\documentclass[11pt,fullpage]{article}

\newcommand{\ol}{\setlength{\itemsep}{0pt.}\begin{enumerate}}
\newcommand{\eol}{\end{enumerate}\setlength{\itemsep}{-\parsep}}
\newcommand{\ignore}[1]{}
\begin{document}
\begin{center}
{\bf Combinatorics hidden in hyperbolic polynomials and related topics 
 }

\vskip 4pt
{Leonid Gurvits }\\
\vskip 4pt
{\tt gurvits@lanl.gov}\\
\vskip 4pt
Los Alamos National Laboratory, Los Alamos , NM 87545 , USA.
\end{center}




\begin{abstract}
The main topic of this paper is various "hyperbolic" generalizations of the Edmonds-Rado theorem on the rank of intersection of two matroids.
We prove several results in this direction and pose a few questions. We also give generalizations of the Obreschkoff theorem and recent results
of J. Borcea and B. Shapiro.

\end{abstract} 
 
 \newtheorem{ALGORITHM}{Algorithm}[section]
\newenvironment{algorithm}{\begin{ALGORITHM} \hspace{-.85em} {\bf :} 
}%
                        {\end{ALGORITHM}}
\newtheorem{THEOREM}{Theorem}[section]
\newenvironment{theorem}{\begin{THEOREM} \hspace{-.85em} {\bf :} 
}%
                        {\end{THEOREM}}
\newtheorem{LEMMA}[THEOREM]{Lemma}
\newenvironment{lemma}{\begin{LEMMA} \hspace{-.85em} {\bf :} }%
                      {\end{LEMMA}}
\newtheorem{COROLLARY}[THEOREM]{Corollary}
\newenvironment{corollary}{\begin{COROLLARY} \hspace{-.85em} {\bf 
:} }%
                          {\end{COROLLARY}}
\newtheorem{PROPOSITION}[THEOREM]{Proposition}
\newenvironment{proposition}{\begin{PROPOSITION} \hspace{-.85em} 
{\bf :} }%
                            {\end{PROPOSITION}}
\newtheorem{DEFINITION}[THEOREM]{Definition}
\newenvironment{definition}{\begin{DEFINITION} \hspace{-.85em} {\bf 
:} \rm}%
                            {\end{DEFINITION}}
\newtheorem{EXAMPLE}[THEOREM]{Example}
\newenvironment{example}{\begin{EXAMPLE} \hspace{-.85em} {\bf :} 
\rm}%
                            {\end{EXAMPLE}}
\newtheorem{CONJECTURE}[THEOREM]{Conjecture}
\newenvironment{conjecture}{\begin{CONJECTURE} \hspace{-.85em} 
{\bf :} \rm}%
                            {\end{CONJECTURE}}
\newtheorem{PROBLEM}[THEOREM]{Problem}
\newenvironment{problem}{\begin{PROBLEM} \hspace{-.85em} {\bf :} 
\rm}%
                            {\end{PROBLEM}}
\newtheorem{QUESTION}[THEOREM]{Question}
\newenvironment{question}{\begin{QUESTION} \hspace{-.85em} {\bf :} 
\rm}%
                            {\end{QUESTION}}
\newtheorem{REMARK}[THEOREM]{Remark}
\newenvironment{remark}{\begin{REMARK} \hspace{-.85em} {\bf :} 
\rm}%
                            {\end{REMARK}}
\newcommand{\alg}{\begin{algorithm}} 
\newcommand{\thm}{\begin{theorem}}
\newcommand{\lem}{\begin{lemma}}
\newcommand{\pro}{\begin{proposition}}
\newcommand{\dfn}{\begin{definition}}
\newcommand{\rem}{\begin{remark}}
\newcommand{\xam}{\begin{example}}
\newcommand{\cnj}{\begin{conjecture}}
\newcommand{\prb}{\begin{problem}}
\newcommand{\que}{\begin{question}}
\newcommand{\cor}{\begin{corollary}}
\newcommand{\prf}{\noindent{\bf Proof:} }
\newcommand{\ethm}{\end{theorem}}
\newcommand{\elem}{\end{lemma}}
\newcommand{\epro}{\end{proposition}}
\newcommand{\edfn}{\bbox\end{definition}}
\newcommand{\erem}{\bbox\end{remark}}
\newcommand{\exam}{\bbox\end{example}}
\newcommand{\ealg}{\end{algorithm}}
\newcommand{\ecnj}{\bbox\end{conjecture}}
\newcommand{\eprb}{\bbox\end{problem}}
\newcommand{\eque}{\bbox\end{question}}
\newcommand{\ecor}{\end{corollary}}
\newcommand{\eprf}{\bbox}
\newcommand{\beqn}{\begin{equation}}
\newcommand{\eeqn}{\end{equation}}
\newcommand{\wbox}{\mbox{$\sqcap$\llap{$\sqcup$}}}
\newcommand{\bbox}{\vrule height7pt width4pt depth1pt}
\newcommand{\qed}{\bbox}
\def\sup{^}
\def\Tp{Tchebyshef polynomial}
\def\Tps{TchebysDeto be the maximafine $A(n,d)$ l size of a code with distance 
$d$hef polynomials}
\newcommand{\rarrow}{\rightarrow}
\newcommand{\larrow}{\leftarrow}
\newcommand{\grad}{\bigtriangledown}

\overfullrule=0pt
\def\setof#1{\lbrace #1 \rbrace}
\section{Introduction, Spectral inequalities}
Consider a homogeneous polynomial $p(x_1,...,x_m)$ of degree $n$ in $m$ real variables.
Such a $p$ is called hyperbolic in the direction $e$ (or $e$- hyperbolic) if for any $x \in R^m$ the polynomial 
 $p(x - \lambda e)$  in the one variable $\lambda$ has exactly $n$ real roots counting their multiplicities. We will assume below
that $p(e) = 1$. Denote an ordered vector of roots of $p(x - \lambda e)$ as 
$\lambda(x) = (\lambda_{1}(x) \geq \lambda_{2}(x) \geq ... \lambda_{n}(x)) $. It is well known that the product of roots is equal
to $p(x)$. Call $x \in R^m$ $e$-positive ($e$-nonnegative) if $\lambda_{n}(x) > 0$ ($\lambda_{n}(x) \geq 0$).
The fundamental result \cite{gar} in the theory of hyperbolic polynomials states that the set of $e$-nonnegative vectors is a 
closed convex cone. A $k$-tuple of vectors $(x_1,...x_k)$ is called $e$-positive ($e$-nonnegative) if $x_i , 1 \leq i \leq k$ are
$e$-positive ($e$-nonnegative). \\

Below, we denote the closed convex cone of $e$-nonnegative vectors as $N_{e}$,
and the open convex cone of $e$-positive vectors
as $C_{e}$. It has been shown in \cite{gar} (see also \cite{khov}) that 
an $e$- hyperbolic polynomial $p$ is also
$d$- hyperbolic for all $e$-positive vectors $d \in C_{e}$. \\

Let us fix $n$ real vectors $x_i \in R^m , 1 \leq i \leq n$ and define the following homogeneous polynomial:
\beqn
P_{x_1,..,x_n}(\alpha_1,...,\alpha_n) = p(\sum_{1 \leq i \leq n} \alpha_i x_i)
\eeqn
Following \cite{khov} , we define the $p$-mixed value of an $n$-vector tuple ${\bf X} = (x_1,..,x_n)$ as
\beqn
M_{p}({\bf X}) = : M_{p}(x_1,..,x_n) = \frac{\partial^n}{\partial \alpha_1...\partial \alpha_n} p(\sum_{1 \leq i \leq n} \alpha_i x_i)
\eeqn
Equivalently, the $p$-mixed value $M_{p}(x_1,..,x_n)$ can be defined by the polarization (see \cite{khov}) :
\beqn
M_{p}(x_1,..,x_n) = 2^{-n} \sum_{b_{i} \in \{-1, +1 \}, 1 \leq i \leq n }  p(\sum_{1 \leq i \leq n } b_i x_i) \prod_{1 \leq i \leq n } b_i
\eeqn

Let us denote as $I_{k,n}$ the set of vectors $r = (r_1,...,r_k)$   with nonnegative integer components and  $\sum_{1 \leq i \leq k} r_i = n$. \\
Let us associate with any such vector $r$ an $n$-tuple of $m$-dimensional vectors  ${\bf X}_{r}$
of  $r_i$ copies of $x_{i}  (1 \leq i \leq k) $.  Notice that 
$$
{\bf X}_{r} = (y_1,...,y_n) ;  y_i \in \{x_1,...,x_k\},  1 \leq i \leq k.
$$
It follows, for instance from the polarization identity (3), that
\beqn
P_{x_1,..,x_n}(\alpha_1,...,\alpha_n) = \sum_{r \in I_{n,n}} \prod_{1 \leq i \leq n } \alpha_{i}^{r_{i}} M_{p}({\bf X}_{r})
\frac{1}{\prod_{1 \leq i \leq n } r_{i}!}
\eeqn
For nonnegative tuple ${\bf X} = (x_1,..,x_n)$, define its capacity as:
\beqn
Cap({\bf X}) = \inf_{\alpha_i > 0, \prod_{1 \leq i \leq n }\alpha _i = 1} P_{x_1,..,x_n}(\alpha_1,...,\alpha_n)
\eeqn
\xam
Probably the best known example of a hyperbolic polynomial is 
\beqn
P(\alpha_0 ,...,\alpha_k) = Det(\sum_{0 \leq i \leq k } \alpha_i A_i) 
\eeqn
where $A_i, 0 \leq i \leq k$ are hermitian matrices and the linear space spanned by $A_i, 0 \leq i \leq k$ contains
a strictly positive definite matrix: $\sum_{0 \leq i \leq k } \beta_i A_i = B \succ 0 $.
This polynomial is hyperbolic in the direction $\beta = (\beta_1,...,\beta_k)$. We can assume wlog
that $B=I$ and that $\beta = (1,0,0,...,0)$. In other words, after  a 
nonsingular linear change of variables
\beqn
P(\alpha_0 ,...,\alpha_k) = Det(\sum_{0 \leq i \leq k } \alpha_{i} B_i) 
\eeqn
where the matrices $B_i, 1 \leq i \leq k$ are hermitian and $B_0 = I$. \\

In this case mixed forms are called mixed discriminants. Let $A_1...A_n$ be $n\times n$  matrices.
The number
\beqn
D(A_1...A_n) = \frac{\partial^n}{\partial x_1...\partial
x_n} \det(x_1 A_1 + \cdots
x_n A_n)
\eeqn
is called the {\it mixed discriminant} of $A_1...A_n$.
\exam
Whether or not
the cone $N_{e}$  of $e$-nonnegative vectors allows a semidefinite representation is a major open question in the area.
In the case of three variables this question was recently settled in \cite{lax} which is a rather direct application of \cite{hel}. 
In this three variables case if
$e = (1,0,0)$ then any $e$- hyperbolic polynomial has a determinantal representation (7), in fact, it even has one with real symmetric
matrices $B_i$. 

\pro
Consider  a homogeneous polynomial $p(x_1,...,x_m)$ of degree $n$ in $m$ real variables which is
hyperbolic in the direction $e$. 
For any pair of $m$-dimensional vectors $x,y \in R^{m}$ there exist a pair of $n \times n$ real symmetric matrices $A,B$
such that $ \lambda(ax + by) = \lambda(aA + bB) $, where $a,b \in R$ and $\lambda(M)$ is the ordered vector of eigenvalues of 
the matrix $M$.
\epro 

\prf
Consider the following  hyperbolic in the direction $(1,0,0)$ polynomial $Q(x_1 ,x_2 ,x_3) = p(x_1 e + x_2 x + x_3 y)$.
Then there exists two $n \times n$ real symmetric matrices $A,B$ such that $p(x_1 e + x_2 x + x_3 y) = Det(x_1 I + x_2 A + x_3 B)$.
It follows that for real $a,b$ the roots of $ax +by$ coincide with eigenvalues of $aA +bB$ and our proof follows directly.
\eprf

Propositin 1.2 allows to "transfer" many "spectral" statements, known for real symmetric matrices, to the context of general
hyperbolic polynomials. (We assume that conditions of Proposition 1.2 are also satisfied in Corollary 1.3 below).
\cor
\begin{enumerate}
\item
Consider a symmetric (i.e. invariant respect to all permutations of variables)$f(y_1,..,y_n): X \rightarrow  R$ ,
where $X \subset R^{n}$ and either $X = R^{n}$, either $X =R^{n}_{+}$ (nonnegative orthant) or $X =R^{n}_{++}$ (positive orthant).
Define ${\bf f}(x) = f(\lambda_{1}(x),\lambda_{2}(x) ... \lambda_{n}(x))$ where either $x \in R^{n}$, either $x \in N_{e}$ or
$x \in C_{e}$ correspondingly.
If $f$ is convex on $X$ then ${\bf f}(x)$ is also convex on either $R^{n}$,  either $N_{e}$ or $C_{e}$ correspondingly. 
(Most recent "hyperbolic inequalities" papers (\cite{kry}, \cite{gul}, \cite{lewis} etc.) 
are simple corollaries of this statement.)
\item
Recall that for two $n \times n$ hermitian $A,B$ there is 
a complete polytope description (Horn's inequalities) \cite{klyach} of
all possible triplets of vectors $(\lambda(A+B),\lambda(A), \lambda(B))$ such that:
$$
\sum_{i \in T}\lambda_{i}(A+B) \leq \sum_{j \in S}\lambda_{j}(A) + \sum_{k \in U}\lambda_{k}(B) ,
$$
where $T,S,U$ are some suitable subsets of $\{1,...,n\}$ of the same size. 

We get from Proposition 1.2 that for any two vectors $x,y \in R^{m}$ the ordered vectors of roots $ 
(\lambda(x+y),\lambda(x), \lambda(y))$ 
satisfy all Horn's inequalities. In particular they satisfy the Lidskii property :
there exists a doubly stochastic $n \times n$ matrix $D$ such that
$\lambda(x+y)-\lambda(x) = D \lambda(y)$. (This settles Open Problem 3.6 posed in \cite{lewis}).
\end{enumerate}
\ecor

\rem In the very same way one can obtain results on
self-concordance. Indeed, results on self-concordance are of the
following nature: consider some symmetric function $f$, and for a pair
$x,y \in R^m$ define $F(t) = {\bf f}(x + th)$.  Assume that $x$
belongs to some cone, usually it is a cone of positive vectors
\cite{gul}.  One needs that F is convex and satisfies some
inequalities for derivatives of F at zero:
$$
|F^{(3)}| \leq 2(F^{(2)})^{1.5};~|F^{(1)}| \leq \sqrt{\sigma} (F^{(2)})^{.5}
$$
Again, if these inequalities and convexity hold for hermitian $n \times n$ matrices, then
we get the same stuff for vectors/$e$-positive vectors/ $e$-positive vectors satisfying $p(x) \geq a > 0$ using
the hyperbolic polynomial $p(x_1 e + x_2 x + x_3 h) = Det(x_1 I + x_2 A + x_3 B)$. \\
The class of inequalities which follow from Proposition 1.2 is larger then the class of symmetric convex inequalities
from \cite{lewis}. For a complex matrix $A$, consider its singular values $a_1 \geq a_2 \geq...\geq a_n \geq 0 $.
Define $f(A) = \sqrt{\sum_{i} i a_i}$. Then $f(A+B) \leq f(A) +f(B)$ \cite{nonhol}, \cite{cdc93}.\\ 
Other inequalities of this type, which are obtained using optimal nonholonomic control, can be found in \cite{cdc93}.
And all of them are non-convex, all of them can be "transfered" to hyperbolic polynomials. \\
Many other things become much more transparent using the polynomial (in three real variables) $p(x_1 e + x_2 x + x_3 h)$.
For instance, the mentioned above fact that $e$- hyperbolic polynomial $p$ is also
$d$- hyperbolic for all $e$-positive vectors $d \in C_{e}$ and the equalities \\
$$
N_{e}=N_{d}, C_{e} = C_{d}
$$
We will get more applications of this ``trick'' (i.e. using hyperbolic polynomial $p(x_1 e + x_2 x + x_3 y)$ in three variables)
in {\bf Section 3}.
\erem

\section {Inequalities for mixed forms, Combinatorics of mixed forms}
In this section we will try to understand if another important class of inequalities, which is valid for determinantal polynomials (7) ,
can be ``transfered" to general hyperbolic polynomials.
We recall below some facts about mixed discriminants, which are mixed forms corresponding to the determinantal polynomials.\\
\begin{description}
\item [ Alexandrov-Fenchel inequalities.] 
Consider an $n$-tuple of positive-semidefinite $n \times n$ hermitian matrices $(A_1,..,A_n)$.
Then the mixed discriminant satisfies the following \\
(hyperbolic) inequality :
\beqn
D(A_1,A_2,A_3,..,A_n)\geq \sqrt{D(A_1,A_1,A_3,..,A_n) D(A_1,A_2,A_3,..,A_n)}
\eeqn
This inequality holds also for mixed forms $M_{p}(x_1,..,x_n)$, where $p$ is $e$-hyperbolic polynomial of degree $n$ ,
and $(x_1,..,x_n)$ are $e$-nonnegative vectors \cite{khov}.
\item [van der Waerden inequality]
Call an $n$-tuple of positive-semidefinite $n \times n$ hermitian matrices $(A_1,..,A_n)$ doubly stochastic if
\beqn
tr(A_i) = 1,~1 \leq i \leq n ; \sum_{1 \leq i \leq n} A_i = I
\eeqn
Then 
\beqn
D(A_1,A_2,A_3,..,A_n) \geq D(\frac{1}{n}I,\frac{1}{n}I ,\frac{1}{n}I ,...,\frac{1}{n}I) = \frac{n!}{n^{n}}
\eeqn
Moreover the inequality is strict if $(A_1,A_2,A_3,..,A_n) \neq (\frac{1}{n}I,\frac{1}{n}I ,\frac{1}{n}I ,...,\frac{1}{n}I) $.\\
The inequality (11) was conjectured in \cite{bapat} and was proved for the
real case in \cite{GS}; the hermitian case and uniqueness were 
proved in \cite{gur}. 
\item [Scaling] Here is the version of (11) which does not require
doubly stochasticity. Similarly to (5) define 
\beqn
Cap(A_1,A_2,A_3,..,A_n) = \inf_{\alpha_i > 0, \prod_{1 \leq i \leq n
}\alpha _i = 1} Det(\alpha_1 A_1,...,\alpha_n A_n) \eeqn Then \beqn
D(A_1,A_2,A_3,..,A_n) \leq Cap(A_1,A_2,A_3,..,A_n) \leq
\frac{n^{n}}{n!} D(A_1,A_2,A_3,..,A_n) \eeqn 
The proof of this
inequality in \cite{GS}, \cite{GS1} requires a matrix scaling. Left
inequality in (13) holds also for mixed forms of nonnegative vectors
for general hyperbolic polynomials.
\item ["Concavity" of $ln(Cap)$ ]
We present below a general result, i.e. which holds for general hyperbolic polynomials.
\lem
Consider an $e$-nonnegative tuple ${\bf X} = (x_1,..,x_n)$. For a vector $r =(r_1,..., r_n) \in I(n,n)$ define
$f(r) = ln(Cap({\bf X}_{r}) $. The function $f(r)$ is concave on $I(n,n)$. I.e.
if $r_0 = \sum_{1 \leq i \leq k} a_i r_i$, where $\sum_{1 \leq i \leq k} a_i = 1 ; a_i \geq 0, 1 \leq i \leq k $ and
$r_i \in I(n,n),  0 \leq i \leq k $, then the following inequality holds :
\beqn
Cap({\bf X}_{r_{0}}) \geq \prod_{1 \leq i \leq k} Cap({\bf X}_{r_{i}})^{a_{i}}
\eeqn
\elem
\prf
We will use a known recent result \cite{kry}, \cite{gul}, \cite{lewis} that the functional $ln(p(x))$ is concave on a positive cone $C_{e}$
(see also Proposition 1.2 and {\bf Section 4 }).
Fix a vector $(\alpha_1,...,\alpha_n) :  \alpha_i > 0, \prod_{1 \leq i \leq n }\alpha _i = 1 $. 

First, let us consider some $z =(z_1,...,z_n) \in I(n,n)$. An easy application of 
the arithmetic/geometric mean inequality
gives that 
\beqn
Cap({\bf X}_{z}) \geq d \mbox{ iff } p(\sum_{1 \leq i \leq n } \alpha_i z_i x_i) \geq d \prod_{1 \leq i \leq n} \alpha_{i}^{z_{i}}
\eeqn
Let $r_i = (r_{i,1},r_{i,2},...,r_{i,n}) ; 0 \leq i \leq k $. It follows from (15) that
$$
ln(p(\sum_{1 \leq j \leq n } \alpha_j r_{i,j} x_j)) \geq \sum_{1 \leq j \leq n } ln(\alpha_j)r_{i,j} + ln(Cap({\bf X}_{r_{i}})), 1 \leq i \leq k
$$
Multiplying the $i$th inequality by $a_i$ and adding afterward we get that
$$
\sum_{1 \leq i \leq k} a_i ln(p(\sum_{1 \leq j \leq n } \alpha_j r_{i,j} x_j)) \geq \sum_{1 \leq j \leq n } ln(\alpha_j)r_{0,j} + 
\sum_{1 \leq i \leq k} a_i ln(Cap({\bf X}_{r_{i}})
$$
Using the concavity of $ln(p(.))$ and (15) we finally get that
$$
Cap({\bf X}_{r_{0}}) \geq \prod_{1 \leq i \leq k} Cap({\bf X}_{r_{i}})^{a_{i}}
$$
\eprf
\item [Edmonds-Rado theorem and Newton polytopes ]
The following result is a direct corollary of the famous Edmonds-Rado theorem on the rank of intersection
of two matroids \cite{gr:lo:sc} :\\
\vskip 5pt
Consider an $n$-tuple of positive-semidefinite $n \times n$ hermitian matrices $(A_1,..,A_n)$.
Then the the mixed discriminant $D(A_1,A_2,A_3,..,A_n) > 0$ iff $Rank(\sum_{i \in S} A_i) \geq |S|$ for
all $S \subset \{1,2,...,n \} $. 

Consider a polynomial $Q(\alpha_1 ,...,\alpha_n) = (Det(\alpha_1 A_1 +...+ \alpha_n A_n)$.
Let us denote as $I_{k,n}$ the set of vectors $r = (r_1,...,r_k)$   with nonnegative integer components and  $\sum_{1 \leq i \leq k} r_i = n$. \\
Associate with an integer vector $r \in I(n,n)$ an $n$-tuple of matrices  ${\bf A}_{r}$ consisting 
of  $r_i$ copies of $A_{i}  (1 \leq i \leq n) $ and denote by $D({\bf A}_{r})$ the corresponding mixed discriminant. Then
\beqn
Q(\alpha_1,...,\alpha_n) = \sum_{r \in I_{n,n}} \prod_{1 \leq i \leq n } \alpha_{i}^{r_{i}} D({\bf X}_{r})
\frac{1}{\prod_{1 \leq i \leq n } r_{i}!}
\eeqn
Therefore the support of $Q$, $supp(Q) = \{ r \in I_{n,n} : D({\bf X}_{r}) > 0 \} $. \\
It follows from Edmonds-Rado theorem that
\beqn
CO(supp(Q)) \cap I(n,n) = supp(Q) ,
\eeqn 
where $CO(supp(Q))$ is a convex hull of $supp(Q)$ , i.e. the Newton polytope of the polynomial $Q$ . \\
The inequality (13) and Lemma (2.1) give a more precise statement :
If $r_0 = \sum_{1 \leq i \leq k} a_i r_i$, where $\sum_{1 \leq i \leq k} a_i = 1, ~a_i \geq 0, ~1 \leq i \leq k $ and
$r_i \in I(n,n),  ~0 \leq i \leq k $, then the following inequality holds:
\beqn
D({\bf X}_{r_{0}}) \geq \prod_{1 \leq i \leq k} D({\bf X}_{r_{i}})^{a_{i}} \frac{n!}{n^{n}}
\eeqn
It is interesting to notice that proofs of  as (11) as well (13) do not use Edmonds-Rado theorem.
\end{description}

Now we are ready to ask relevant questions for hyperbolic polynomials.
\begin{description}
\item[Question 1.]
Consider  a homogeneous polynomial $p(y_1,...,y_m)$ of degree $n$ in $m$ real variables which is
hyperbolic in the direction $e$ and the corresponding polynomial in $n$ real variables
$$
P_{x_1,..,x_n}(\alpha_1,...,\alpha_n) = p(\sum_{1 \leq i \leq n} \alpha_i x_i)
$$
where $x_1,..,x_n \in R^m$ are $e$-nonnegative. Is it true that
$$
CO(supp(P_{x_1,..,x_n})) \cap I(n,n) = supp(Q) ?
$$
(Recall that the convex hull $CO(supp(P_{x_1,..,x_n}))$ 
is the Newton polytope of the polynomial $P_{x_1,..,x_n}$.)\\
If the answer is ``yes" then we get an analog of Edmonds-Rado theorem for hyperbolic polynomials.
To state this, 
define the $p$-rank of $x \in R^m$ as $Rank(x) = |\{i : \lambda_{i}(x) \neq 0 \} |$.
Then, the statement is that
for  $e$-nonnegative tuples ${\bf X} = (x_1,..,x_n)$ the $p$-mixed value $M_{p}(x_1,..,x_n) > 0$ 
iff $Cap({\bf X}) > 0$; or equivalently, iff  $Rank(\sum_{i \in S} x_i) \geq |S|$ 
for all $S \subset \{1,2,...,n \} $. 

This (conditional) result follows from the following {\bf Proposition 2.2} and {\bf Proposition 2.3} .
\pro
Consider  a homogeneous $e$-hyperbolic polynomial $p(.)$ of degree $n$ in $m$ real variables.
Let ${\bf X} = (x_1,..,x_n)$ be an $e$-nonnegative tuple.
Then $Cap({\bf X}) > 0$ iff the following generalized Edmonds-Rado condition holds:\\
$Rank(\sum_{i \in S} x_i) \geq |S|$ for all $S \subset \{1,2,...,n \} $.
\epro
\prf
We will use two known facts: \\
Fact 1. $\lambda_{k}(x+y) \geq \lambda_{k}(x) $ provided $y$ is 
$e$-nonnegative; if $x$ is $e$-nonnegative and the scalar $a > 0$ then
$\lambda_{k}(a x) = a \lambda_{k}(x)$, $ 1 \leq k \leq n$.\\
Fact 2. $p(x) = p(e) \prod_{1 \leq k \leq n} \lambda_{k}(x)$; if $x, y$ are $e$-nonnegative then $p(x+y) \geq p(x)$. \\

If the generalized Edmonds-Rado condition holds then for any subset $S$ of cardinality $k$ we have the inequality
$$
\lambda_{k}(\sum_{i \in S} x_i) \geq r_k > 0
$$
Using Fact 1 and the same argument as in the proof of Theorem 1.14 from \cite{GS1}, we get that for $e$-positive  ($z_i, ~1 \leq i \leq n $)
the following inequality holds:
$$
p(\sum_{1 \leq i \leq n} z_i x_i) \geq \prod_{1 \leq k \leq n} r_{k} z_{k}
$$
In other words, $Cap({\bf X}) \geq \prod_{1 \leq k \leq n} r_{k} $. \\
Suppose that  the generalized Edmonds-Rado condition does not hold ; or, wlog ,
$\lambda_{k}(\sum_{1 \leq i \leq k} x_i) = 0$. Also, assume wlog that all vectors $e-x_i$ are $e$-nonnegative ,
where $e$ is a vector (direction) used in the definition of hyperbolicity (recall that $p(e)=1$). Choose $z_i = a > 0$ for $1 \leq i \leq k$
and $z_i = b > 0$ for $k+1 \leq i \leq n$. Using Fact 2, we get the following inequality:
$$
p(\sum_{1 \leq i \leq n} z_i x_i) \leq p(\sum_{1 \leq i \leq k} a x_i + (n-k) b e)
$$
By our assumption, an $e$-nonnegative vector $\sum_{1 \leq i \leq n}  x_i$ has at most $k-1$ nonzero roots: \\
$ r \geq r_1 \geq...\geq r_{k-1} \geq 0$. Therefore
$$
p(\sum_{1 \leq i \leq k} a x_i + (n-k) b e) =  b^{n-k+1} \prod_{1 \leq i \leq k-1} a r_i + b \leq  b^{n-k+1} (ar + b)^{k-1}
$$
Finally,
$$
\frac{p(\sum_{1 \leq i \leq n} z_i x_i)}{\prod_{1 \leq i \leq n} z_i } \leq \frac{(ar + b)^{k-1} b}{a^{k}}
$$
For a fixed $a > 0$ the right side of the last inequality converges to zero if (positive) b converges to zero.
Therefore, $Cap({\bf X}) = 0$.
\eprf

\pro
For any vector $r =(r_1,..., r_n) \in I(n,n)$ the capacity
$Cap({\bf X}_{r}) > 0 $ iff $r \in CO(supp(P_{x_1,..,x_n}))$.
\epro
\prf
It is an easy application of convexity of the logarithm on the positive orthant and the Hahn-Banach separation theorem; all what is important
is that the coefficients of $P_{x_1,..,x_n}$ are nonnegative.
\eprf
\dfn
Consider a polynomial $p(y_1,...,y_m)$ of degree $n$ in $m$ real variables
hyperbolic in direction $e$ and assume
that $P(e)=1$. Call an $n$ -tuple ${\bf X} = (x_1,..,x_n)$ of $m$-dimensional real vectors $d$-doubly stochastic if
$x_i, 1 \leq i \leq n$ are $e$-nonnegative, $\sum_{1 \leq i \leq n} x_i = d \in C_{e}$ and $tr_{d}(x_i) = 1, 1 \leq i \leq n$ ,
where $tr_{d}(x)$ is the sum of roots of $x$ in the direction $d \in C_{e}$, 
i.e. roots of the equation $p(x - td) = 0$.
\edfn

\pro
If an $n$-tuple ${\bf X}$ is $d$-doubly stochastic then $Cap({\bf X})= p(d)$
\epro

\prf
Since $p(\sum_{1 \leq i \leq n} x_i) = p(d)$ hence $Cap({\bf X}) \leq p(d) $.
It remains to prove that $p(\sum_{1 \leq i \leq n} e^{a_{i}} x_i) \geq p(d)$ if $\sum_{1 \leq i \leq n} a_i = 0$
and $a_i$ are real. The functional $g(a_1,...,a_n) = p(\sum_{1 \leq i \leq n} e^{a_{i}} x_i)$ is convex on $R^{n}$
(even $log(g(a_1,...,a_n))$ is convex) as all  coefficients of $P_{x_1,..,x_n}$ are nonnegative.\\
Thus we need to check, similarly to the proof of Lemma 3.6 in \cite{GS1}, that the gradient of $g$ evaluated
at the zero vector $(0,...,0)$ is a constant multiple of $(1,...,1)$. But at the zero vector,
$$
\frac{\partial}{\partial a_k } (g(a_1,...,a_n)) = \frac{\partial}{\partial a_i } p (\sum_{1 \leq i \leq n} (1 + a_{i}) x_i), 1 \leq k \leq n
$$
Using  $d$-double stochasticity, we get that 
$$
p (\sum_{1 \leq i \leq n} (1 + a_{i}) x_i) = p(d + \sum_{1 \leq i \leq n}  a_{i} x_i)
$$
Finally, at the zero vector, we get that
$$
\frac{\partial}{\partial a_k } (g(a_1,...,a_n)) = \frac{\partial}{\partial a_k } p(d + a_{k} x_k) = p(d)^{-1} tr_{d}(x_k)= p(d)^{-1}, 1 \leq k \leq n
$$
Therefore, the zero vector $(0,...,0)$ is a global (not always unique) minimum of $g(a_1,...,a_n)$ on the hyperplane
$ (a_1,..., a_n): \sum_{1 \leq i \leq n} a_i = 0 $.  This means that
$$
p(\sum_{1 \leq i \leq n} \alpha_i x_i) \geq p(d)\prod_{1 \leq i \leq n}  \alpha_i  ; \alpha_i > 0, 1 \leq i \leq n
$$
Thus $Cap({\bf X})= p(d)$.
\eprf

\rem
Perhaps, Proposition 2.5 can be strengthened to the following statement: \\
let $\Lambda = (\lambda_1,... ,\lambda_n)$ be roots of $\sum_{1 \leq i \leq n} c_i x_i$ in the direction $d \in C_{e}$, where 
the tuple $(x_1,...,x_n)$ is $d$-doubly stochastic and $(c_i ,1 \leq i \leq n)$ are real numbers.
Then there exists a doubly stochastic $n \times n$ matrix $D$ such that $\Lambda = C D$, where
$C= (c_1,..., c_n)$.\\
\erem

\item[Question 2 ] 
Define the van der Waerden constant of a hyperbolic polynomial $p(y_1,...,y_m)$ of degree $n$ in $m$ real variables as
$$
VDW(p) = \inf \frac{M_{p}(x_1,..,x_n)}{Cap(x_1,..,x_n)}
$$
where the infimum is taken over the set of tuples $(x_1,..,x_n)$ of
$e$-positive vectors. It is easy to see that $VDW(p) \leq
\frac{n!}{n^{n}}$. Is $VDW(p) = \frac{n!}{n^{n}}$ ? Is it positive ?
\\ This question is a ``hyperbolic" analog of the van der Waerden
conjecture for permanents/mixed discriminants.  If the van der Waerden
constant is positive then our analog of the Edmonds-Rado theorem for
hyperbolic polynomials follows. \\
\end{description}

\subsection{``Hyperbolic" scaling }
Let us explain why {\bf Question 2} above is indeed an analog of the van der Waerden conjecture for hyperbolic polynomials.
\lem
The infimum in (5) is attained iff there exist positive numbers $(\alpha_1,...,\alpha_n)$ with $\prod_{1 \leq i \leq n}\alpha_i = 1$ 
 and an $e$-positive vector $d \in C_{e}$
such that the tuple $(\alpha_1 x_1,...,\alpha_n x_n)$ is $d$-doubly stochastic.
\elem

\prf 
The ``if" part follows directly from Proposition 2.5. Moreover in this case $Cap({\bf X}) = p(d)$.
The ``only if" part follows, very similarly to the proof of Proposition 2.5, from the standard
necessary condition for the corresponding conditional extremum.
\eprf

A fairly direct adaption of Lemmas 3.7, 3.8 from \cite{GS1} gives that The infimum in (5) is attained and unique
if a tuple ${\bf X}$ is $e$-positive. Therefore Question 2 is equivalent to finding
$$
\inf \frac{M_{p}(x_1,..,x_n)}{p(d)}: (x_1,..,x_n) \mbox{ is } d-\mbox{doubly stochastic}, d \in C_{e}
$$
\dfn
Consider an $e$-nonnegative tuple ${\bf X} = (x_1,..,x_n)$ such that the sum of its components
 $S({\bf X}) = d = \sum_{1 \leq i \leq k} x_i$ is $e$-positive.
Define the following map (Hyperbolic Sinkhorn) acting on such tuples:
$$
HS({\bf X}) = {\bf Y} = (\frac{x_{1}}{tr_{d}(x_{1})},..., \frac{x_{n}}{tr_{d}(x_{n})})
$$
Hyperbolic Sinkhorn Iteration ({\bf HSI}) is a recursive procedure:
$$
{\bf X}_{j+1} = HS({\bf X}_{j}), j \geq 0, ~{\bf X}_{0} \mbox{ is an $e$-nonnegative tuple with } \sum_{1 \leq i \leq k} x_i \in C_{e}\;.
$$
We also define the doubly-stochastic defect of $e$-nonnegative tuples with $e$-positive sums as
$$
DS({\bf X}) = \sum_{1 \leq i \leq k} (tr_{d}(x_{i}) - 1)^{2} ; \sum_{1 \leq i \leq k} x_i = d \in C_{e}
$$
\edfn
\xam
Consider the following hyperbolic polynomial in $n$ variables: $p(z_1,...,z_n) = \prod_{1 \leq i \leq n} z_i$.
It is $e$- hyperbolic for $e = (1,1,...,1)$. And $N_{e}$ is a nonnegative orthant, $C_{e}$ is a positive orthant.
An $e$-nonnegative tuple ${\bf X} = (x_1,..,x_n)$ can be represented by an $n \times n$ matrix $A_{{\bf X}}$ with nonnegative entries:
the $i$th column of $A$ is a vector $x_i \in R^{n}$. If $Z = (z_1,...,z_n) \in R^{n}$ and
$d = (d_1,...,d_n) \in R^{n} ; z_i > 0, 1 \leq i \leq n$, then $tr_{d}(Z) = \sum_{1 \leq i \leq n} \frac{z_{i}}{d_{i}}$.

Recall that for a square matrix $A = \{a_{ij}: 1 \leq i,j \leq N\}$ row
scaling is
defined as
$$
R(A) = \{ \frac{a_{ij}}{\sum_j a_{ij}} \}, $$
column scaling as $C(A) = \{ \frac{a_{ij}}{\sum_i a_{ij}} \}$
assuming that all
denominators are nonzero.
The iterative process $...CRCR(A)$ is called {\em Sinkhorn's iterative scaling} (SI).
In terms of the matrix $A_{{\bf X}}$ the map $HS({\bf X})$ can be realized as follows:
$$
A_{HS({\bf X})} = C(R(A_{{\bf X}}))
$$
So, the map $HS({\bf X})$  is indeed a (rather far-reaching) generalization of Sinkhorn's scaling. Other generalizations (not all hyperbolic)
can be found in \cite{GY}, \cite{stoc}, \cite{arxiv}.
\exam
Before proving the next theorem let us first state and prove some properties of the map $HS$.
\lem
Consider  an $e$-nonnegative tuple ${\bf X} = (x_1,..,x_n)$ with $S({\bf X}) \in C_{e}$, i.e. $S({\bf X})$ being $e$-positive.
Then $Cap(HS({\bf X})) = (\prod_{1 \leq i \leq n} tr_{S({\bf X})}(x_{i}))^{-1} Cap({\bf X}) $, 
and $p(S(HS({\bf X}))) \leq p(S({\bf X}))$.
\elem

\prf
Consider the following polynomial in one variable $D(t) = p(td+x) = \sum_{0 \leq i \leq n} c_{i} t ^{i} $.
It follows from the identity (4) that \\
\beqn
c_{n} = M_{p}(d,..,d) (n!)^{-1} = p(d), c_{n-1} = M_{p}(x,d,..,d) (1! (n-1)!)^{-1},..., c_{0} = M_{p}(x,..,x) (n!)^{-1} = p(x).
\eeqn
Let $(\lambda_{1}^{(d)}(x) \geq \lambda_{2}^{(d)}(x) \geq... \geq \lambda_{n}^{(d)}(x))$ be the (real) roots of $x$
in the $e$-positive direction $d$, i.e. the roots of the equation $p(td-x) = 0$. 
Notice that
$$
\prod_{1 \leq i \leq n} \lambda_{i}^{(d)}(x) = \frac{c_{0}}{c_{n}} \mbox{ and } tr_{d}(x) =\frac{c_{n-1}}{c_{n}}.
$$

Thus we get the following identities:
\beqn
p(x) = p(d) \prod_{1 \leq i \leq n} \lambda_{i}^{(d)}(x), tr_{d}(x) =  M_{p}(x,d,..,d) ((n-1)! p(d))^{-1}
\eeqn
The first statement follows directly from the following obvious formula
\beqn
Cap(a_1 x_1 ,..., a_n x_n) = (\prod_{1 \leq i \leq n} a_i ) Cap(x_1,...,x_n) ; a_i > 0\;. 
\eeqn
The second identity in (20) reproves the
well known fact that the functional $tr_{d}(x)$ is linear.
Therefore, we get that 
$$
tr_{S({\bf X}))}(HS({\bf X})) = \sum_{1 \leq i \leq n} tr_{d} (\frac{x_{i}}{tr_{d}(x_{i})}) = n\;.
$$

Via the standard arithmetic/geometric means inequality and using the first identity in (20),
we finally get that
$$
p(S(HS({\bf X}))) = p(d) \prod_{1 \leq i \leq n} \lambda_{i}^{(d)}(S(HS({\bf X}))) \leq 
p(d) (\sum_{1 \leq i \leq n} \lambda_{i}^{(d)}(S(HS({\bf X}))))^{\frac{1}{n} } = p(d)\;.
$$
\eprf

We also need the following ``quantitative" version of Proposition 2.5.
\pro
Consider an $e$-nonnegative tuple ${\bf X} = (x_1,..,x_n)$ with $d = \sum_{1 \leq i \leq n} x_i$ being $e$-positive.
If $DS({\bf X} ) = \sum_{1 \leq i \leq k} (tr_{d}(x_{i}) - 1)^{2} \leq \frac{1}{n} $ then $Cap({\bf X} ) > 0 $
\epro

\prf
Quite naturally in this context (see, for instance, \cite{lsw}), 
we will use Proposition 2.2 , i.e.
we will prove the the conditions of this proposition imply 
the generalized Edmonds-Rado condition : \\
$Rank(\sum_{i \in S} x_i ) \geq |S|$ for all $S \subset \{1,2,...,n \} $. 

Suppose that  the generalized Edmonds-Rado condition does not hold or, wlog,
suppose 
$\lambda_{k}^{(d)}(\sum_{1 \leq i \leq k} x_i) = 0$, where $1 \leq k \leq n-1$.\\
Since $N_{e}=N_{d}, C_{e} = C_{d}$,  we can use Facts 1 and 2 
stated in the proof of Proposition 2.2,
to get that
$$
\lambda_{i}^{(d)}(\sum_{1 \leq i \leq k} x_i) \leq \lambda_{i}^{(d)}(d) = 1, 1 \leq i \leq k-1 
$$
Therefore, $tr_{d}(\sum_{1 \leq i \leq k} x_i) \leq k-1$. On the other hand, using the linearity of 
the functional $tr_{d}(x)$, we obtain that 
$$
tr_{d}(\sum_{1 \leq i \leq k} x_i) = \sum_{1 \leq i \leq k} tr_{d}(x_{i}) = \sum_{1 \leq i \leq k} (1 - \delta_{i}),
$$
where $\sum_{1 \leq i \leq n} (\delta_{i})^{2} \leq \frac{1}{n} $.
Therefore, the Cauchy-Schwarz inequality implies that 
$$
\sum_{1 \leq i \leq k} |\delta_{i}| \leq \sqrt{\frac{k}{n}} < 1 
$$
This gives that $tr_{d}(\sum_{1 \leq i \leq k} x_i) > k-1$, the desired contradiction. \\
\eprf

\thm
Consider Hyperbolic Sinkhorn Iteration ({\bf HSI}) ${\bf X}_{j+1} = HS({\bf X}_{j}); ~j \geq 0 $, where the initial
$e$-nonnegative tuple ${\bf X}_{0} = (x_1,..,x_n)$ satisfies $d_{0} = \sum_{1 \leq i \leq n} x_i \in C_{e}$.
Then the following statements are equivalent:
\begin{enumerate}
\item The exists $j \geq 0$ such that $DS({\bf X}_{j}) \leq \frac{1}{n} $
\item $Cap({\bf X}_{0}) > 0$
\item $\lim DS({\bf X}_{j}) = 0 $
\end{enumerate}
\ethm

\prf
The implication $1 \rightarrow 2$ is Proposition 2.11, the implication $3 \rightarrow 1$ is obvious.
It remains to prove that $2 \rightarrow 3$. Let us introduce the following notations: \\
${\bf X}_{j} = (x_{j,1},...,x_{j,n})$, $\sum_{1 \leq i \leq n} x_{j,i} = d_j$,
$tr_{d_{j}}(x_{j,i}) = a(j,i)$, $\prod_{1 \leq i \leq n} a(j,i) = F_j$. 

First, Lemma 2.10 gives that $p(d_{j+1}) \leq p (d_{j}) \leq p(d_{0})$. Thus, directly from the definition (5),
$Cap({\bf X}_{j}) \leq p(d_{0}) < \infty $. In other words, the sequence $(Cap({\bf X}_{j}), j \geq 0)$ is bounded. \\
By the definition of  Hyperbolic Sinkhorn Iteration ({\bf HSI}) we get that
$$
x_{j+1,i} = \frac{x_{j,i}}{a(j,i)}, \sum_{1 \leq i \leq n} a(j,i) = 1 
$$
Therefore, using (21) and the arithmetic/geometric means inequality, we obtain that
$$
Cap({\bf X}_{j+1}) = F_{j}^{-1} Cap({\bf X}_{j}), F_{j}^{-1} \geq 1 ; j \geq 0
$$
Moreover, if $DS({\bf X}_{j})$ does not converge to zero then the product $ P_{j} = \prod_{0 \leq k \leq j}F_{j}^{-1}$
converges to infinity.$Cap({\bf X}_{j+1}) = P_{j} Cap({\bf X}_{0})$ and $Cap({\bf X}_{0}) > 0$, 
therefore if $DS({\bf X}_{j})$ does not converge to zero the sequence $(Cap({\bf X}_{j}), j \geq 0)$ is not bounded.
This is the desired contradiction.
\eprf

\rem
We can define the map $HS(.)$ directly in terms of the polynomial 
$$
Q(\alpha_1,...,\alpha_n) = P_{x_1,..,x_n}(\alpha_1,...,\alpha_n) = p(\sum_{1 \leq i \leq n} \alpha_i x_i).
$$
Indeed, if $ \sum_{1 \leq i \leq n} \alpha_i x_i =d \in C_{e}$ then
\beqn
tr_{d}(\alpha_i x_i) = \frac{ \alpha_{i } \frac{\partial}{\partial \alpha_{i } } Q(\alpha_1,...,\alpha_n) } {Q(\alpha_1,...,\alpha_n)}
\eeqn
As $Q(\alpha_1,...,\alpha_n)$ is a homogeneous polynomial of degree $n$ thus it satisfies Euler's identity:
$$
Q(\alpha_1,...,\alpha_n) = n \sum_{1 \leq i \leq n} \alpha_i \frac{\partial}{\partial \alpha_{i } } Q(\alpha_1,...,\alpha_n)
$$
(Notice that the linearity of $tr_{d}(x)$ for $d \in C_{e}$ is a particular case of  Euler's identity.)

Using formula (22), we can redefined the map $HS(.)$ as
$$
F(\alpha_1,...,\alpha_n) = (\frac{Q(\alpha_1,...,\alpha_n)}{\frac{\partial}{\partial \alpha_{1} } Q(\alpha_1,...,\alpha_n) },...,
\frac{Q(\alpha_1,...,\alpha_n)}{\frac{\partial}{\partial \alpha_{n} } Q(\alpha_1,...,\alpha_n) }).
$$
Correspondingly, the inequality $p(S(HS({\bf X}))) \leq p(S({\bf X}))$ can be rewritten as
\beqn
Q((\frac{\partial}{\partial \alpha_{1} } Q(\alpha_1,...,\alpha_n))^{-1},..., (\frac{\partial}{\partial \alpha_{n} } Q(\alpha_1,...,\alpha_n))^{-1})
\leq Q(\alpha_1,...,\alpha_n)^{-(n-1) } ; \alpha_i > 0 
\eeqn
where the equality is achieved iff $ \alpha_{i } \frac{\partial}{\partial \alpha_{i} } = Q(\alpha_1,...,\alpha_n)$. \\ 
The inequality (23) is indeed ``hyperbolic", it is not valid for general homogeneous polynomials with nonnegative coefficients.
Consider $Q(\alpha_1,\alpha_2) = \alpha_{1}^{2} +  \alpha_{2}^{2};~ n =2$. Then
$$
Q((\frac{\partial}{\partial \alpha_{1} } Q(\alpha_1,...,\alpha_n))^{-1},..., (\frac{\partial}{\partial \alpha_{n} } Q(\alpha_1,...,\alpha_n))^{-1})
\geq  Q(\alpha_1,...,\alpha_n)^{-1}  
$$
There is another inequality for general homogeneous polynomials with nonnegative coefficients involving partial derivatives:  the
Baum-Snell-Bregman inequality \cite{baum}, \cite{soules}, \cite{breg-minc}. It is interesting that in the case of Example 2.9 
the Baum-Snell-Bregman inequality is equivalent to (23). Also, in the case of Example 2.9 the map $HS(.)$ is a composition 
of two Bregman's projections
associated with one convex functional $f(x_1,...,x_k) = \sum_{1 \leq i \leq k} x_i ln(x_i) ; x_i \geq 0$ (\cite{breg-sink}, \cite{GY}).
It remains to be 
understood whether the map $HS(.)$ for general hyperbolic polynomials has some  Bregman's projections interpretation.
\erem

\section{Obreschkoff theorem and hyperbolic determinantal polynomials in three variables }
Recall that the companion matrix $C_q$ of a monic polynomial $q(x)=x^{n}-a_1 x^{n-1}-..-a_n$ is
a $n \times n$ matrix defined as follows:
$$
C_q = \left(\begin{array}{ccccc}
		  0 &  1 & 0 &\dots & 0\\
		  0 & 0 & 1 & \dots & 0\\
		  0  & 0 & 0 & 1\dots & 0\\
		  \dots  & \dots & \dots & \dots & \dots \\
		  a_n & a_{n-1} & a_{n-2} & \dots & a_1\end{array} \right).
$$
Consider two monic polynomials of degree $n$, $q = x^{n}-a_1 x^{n-1}-..-a_n$ and $r = x^{n}-b_1 x^{n-1}-..-b_n$, and define
the following homogeneous polynomial of degree $n$ in three real variables:
\dfn
$p(x,y,t) = Det(x C_q  + y C_r -t I)$
\edfn
Notice that with respect to this polynomial  the roots of a vector $(x,y,z)$ with $x+y \neq 0$ in the direction $(0,0,1)$ are
$((x+y)\lambda_1 + z, (x+y)\lambda_2 + z,..., (x+y)\lambda_n + z)$, 
where $(\lambda_1 ,\lambda_2 ,..., \lambda_n)$ are the roots of the polynomial $x q + y r$.
\pro
The polynomial $p(x,y,t)$ is $e$-hyperbolic, where $e= (0,0,1)$, iff all polynomials of the form $\{x p + y r: (x,y) \neq 0 ; x,y \in R \}$
have only real roots.
\epro

\prf
First, let us prove the ``if" part.
Recall that the eigenvalues of the companion matrix $C_q$ are exactly the roots of the polynomial $q$.
Therefore, if $x+y \neq 0$ then the eigenvalues of $x C_q  + y C_r$ are (all real) roots of $x p + y r$ multiplied by a real number $x+y$.
If $x+y = 0$ then the eigenvalues of $x C_q  + y C_r$ are $(0,0,...,x a_1 + y b_1)$ ; and thus also real. 

Second, we prove the 
``only if" part. In the view of the first part we need only to prove that all roots of
polynomial $q-r$ of degree $n-1$ are real. Assume, wlog, that $q$ and $r$ don't have common roots.
Suppose that there exists a complex $z_{0} = x + iy, y > 0$ such that $(q-r)(z) =0$.
In other words the rational nonconstant function $R(z) = \frac{q(z)}{r(z)} - 1$ has a zero in the upper half-plane $H_{+} = \{z: Im(z) > 0 \}$.
Since $R(z)$ is analytic and nonconstant on $H_{+}$ and $R(z_{0})= 0$, \
the range $\{ R(z): |z-z_{0}| \leq \epsilon \}$
contains a complex ball $z: |z| \leq \delta > 0$ for all small enough $\epsilon$. Therefore there exists $z_{1}$ with $Im(z_{1}) > 0$
such that $\frac{q(z_{1})}{r(z_{1})} = 1 + \delta $. It follows that the polynomial $q - (1 + \delta) r$ has a non-real root,
but in this case $1 - (1 + \delta) = -\delta \neq 0$. We got the desired contradiction.
\eprf

\cor
Consider the following ``line" of monic polynomials: $P_{a}(x) = a q(x + b + c a) + (1-a)r(x + b + c a)$,
where $a \in R$ and $b,c$ are fixed real numbers. Let $\lambda_{a}(1) \geq... \geq \lambda_{a}(n)$ be the
roots of $P_{a}$. Let $f(x_1,...,x_n)$ be any symmetric and convex on $R^n$ functional.
Define $F(a)= f(\lambda_{a}(1),...,\lambda_{a}(n))$. If all polynomials of the form $\{x p + y r: (x,y) \neq 0 ; x,y \in R \}$
have only real roots then the function $F(.)$ is convex on $R$.
\ecor
\prf
In terms of the $(0,0,1)$-hyperbolic polynomial $p(x,y,t) = Det(x C_q  + y C_r -t I)$ the roots of
the polynomial $P_{a}$ are equal to the roots of the vector $(a, 1-a, b + ca)$ as $a + (1-a) = 1$.
It remains to apply either \cite{lewis} or Proposition 1.2.
\eprf

\xam
Consider an arbitrary monic polynomial $q$ of degree $n$ with all real roots,
define $P_{a}(x) = q(x + b + c a) + aq^{\prime}(x + b + c a)$, where $q^{\prime}$ is a derivative of $q$.
Then $r = q + q^{\prime}$ is also monic, and $P_{a}(x) = (1-a)q(x + b + c a) + ar(x + b + c a)$.
A well known result gives that that the pair $(q,r)$ satisfies the conditions of Corollary 3.3.
Let $\lambda_{a}(1) \geq... \geq \lambda_{a}(n)$ be the
roots of $P_{a}$. We get that $f_{k}(a) = \sum_{1 \leq i \leq k \leq n} \lambda_{a}(i)$ is a convex
function on $R$, $g_{k}(a) = - \sum_{1 \leq i \leq k \leq n} \lambda_{a}(n-i)$ is a convex
function on $R$.
If $q$ is a monic polynomial $q$ of degree $n$ with $n$ distinct real roots and $P_{a}(x)= q(x+a)- aq^{\prime}(x + a)$,
then $P_{a}$ has $n$ disinct real roots for all $a \in R$.
Thus $f_{k}(a)$ is differentiable for all $k \leq n$ ; a direct application
of the formula for the derivative of implicit functions gives that $f_{k}^{\prime}(0) = 0$. Together with the convexity
it gives that the global minimum of $f_{k}(a)$ is attained at zero, $1 \leq k \leq n$. It is easy to see that $f_{n}(a)$, 
which is the sum of the roots of $P_{a}$, is constant on $R$. Therefore, by a well known result, if $ab \geq 0$ and $|a| \leq |b|$
then there exists a doubly stochastic matrix $D_{a,b}$ such that $\Lambda_{a} = D_{a,b}\Lambda_{b}$ (i.e.
the vector $\Lambda_{a}$ is majorized by $\Lambda_{b}$).

The case of nondistinct
roots can be now proved by a standard perturbation argument: if a sequence of functions $f_{m}: R \rightarrow R, m \geq 1$
converges pointwise on $R$ to a function $f: R \rightarrow R$  and $f_{m}(x) \geq f_{m}(0) ; x \in R,  m \geq 1$ 
then the inequality $f(x) \geq f(0) ; x \in R $ also holds. 

The results from this example generalize some results of the recent paper \cite{shapiro} and simplify proofs of others.

Our solution of Open Problem 3.6 posed in \cite{lewis} (see Part 2 of of Corollary 1.3) provides the following general majorization result:
\cor
Consider a hyperbolic pair $(q,r)$ of monic polynomials of degree $n$, i.e.
a pair $(q,r)$ such
that all polynomials of the form $\{x p + y r: (x,y) \neq 0 ; x,y \in R \}$
have only real roots. Consider two real $3$-dimensional vectors $X = (x,y,z)$ and $\Delta = (\delta_{1},\delta_{2},\delta_{3})$.
Assume that $x+y = L \neq 0,~ \delta_{1} + \delta_{2} = M \neq 0,~ x+y + \delta_{1} + \delta_{2}   = K \neq 0$.
Define the following polynomials 
\begin{eqnarray}
P_{X}(t) = xq(t -\frac{z}{L}) + yr(t-\frac{z}{L}) \nonumber \\
 P_{X+\Delta}(t) = (x+ \delta_{1}) q(t - \frac{\delta_{3}}{K}) + (y + \delta_{2}) 
r(t - \frac{\delta_{3}}{K}) \nonumber \\
P_{\Delta}(t) = \delta_{1}q(t - \frac{\delta_{3}}{M}) + \delta_{2} r(t - \frac{\delta_{3}}{M}). \nonumber
\end{eqnarray}
Let $\Lambda_{X}, \Lambda_{X+\Delta}, \Lambda_{\Delta} $ be the ordered vectors (from the largest to the smallest)of roots of the 
degree-$n$ polynomials 
$P_{X}, P_{X+\Delta}, P_{\Delta} $ correspondingly. Define an $n$-dimensional vector $ORD_{X}$ as an ordering of the vector $L \Lambda_{X}$,
$ORD_{X+ \Delta }$ as an ordering of the vector $K \Lambda_{X+ \Delta }$, $ORD_{\Delta }$ as an ordering of the vector $M \Lambda_{\Delta }$. 
Then the vector $ORD_{X+ \Delta } - ORD_{\Delta }$ is majorized by $ORD_{X}$.
\ecor
(Of course, we can now prove many statements of this kind via applying the
Second Part of Corollary 1.3 in its full generality,
i.e. using all Horn's inequalities.)
\exam

Let us consider two polynomials $q,r$ with real coefficients. Assume that $q$ is monic, the 
degree $deg(q)$ of $q$ is $n$
and also that $q$ has $n$ distinct real roots $\lambda_1 >... > \lambda_n $. If $deg(r) \leq n$ then
\beqn
\frac{r(z)}{q(z)} = A + \sum_{1 \leq k \leq n} \frac{a_{k}}{z -\lambda_{k}} ; A, a_k \in R\;.
\eeqn
If a complex number $z= u + v i$ and $Re(z)=u, Im(z) = v \neq 0 $  then 
\beqn
Im(\frac{r(z)}{q(z)}) = - \sum_{1 \leq k \leq n}  \frac{a_{k} v }{ (u -\lambda_k)^{2} + v^{2} }
\eeqn
Call a pair of polynomials $(q,r)$ hyperbolic if all polynomials of the form $\{x p + y r: (x,y) \neq 0 ; x,y \in R \}$
have only real roots. As explained (and used) in the proof of Proposition 3.2, the hyperbolicity
of a pair of polynomials $(q,r)$ is equivalent to the property
\beqn
Im(\frac{r(z)}{q(z)}) \neq 0 \mbox{  if  } Im(z) \neq 0\;.
\eeqn
Therefore if all $a_k$ in () are of the same sign, say $a_{k} \geq 0$, then the pair of polynomials $(q,r)$ is hyperbolic.
Moreover in this case $n \geq deg(r) \geq n-1$ ; if the pair$(q,r)$ is coprime (i.e. no common roots) then
 all polynomials of the form $\{x p + y r: (x,y) \neq 0 ; x,y \in R \}$
have distinct roots  as in this case $0 \neq a_{k} = \frac{r(\lambda_{k})}{q^{\prime}(\lambda_{k})},~ 1 \leq k \leq n $.
(Recall that the condition $ \frac{r(\lambda_{k})}{q^{\prime}(\lambda_{k})} > 0, 1 \leq k \leq n $ forces the
interlacing of the roots.) 

What we proved above is a slightly less general (because of the assumption that the roots of $q$ are distinct) sufficiency
part of the Obreschkoff theorem \cite{obres}. We will prove below the following analytic version of the necessity part.
\thm
Consider two analytic functions $F,G$ with real Taylor series. Assume that all roots of $F$ are real and simple:
i.e. the set of roots of $F$ is $(\lambda_k \in R, -\infty \leq L < k < U \leq \infty)$ and $F^{\prime}(\lambda_{k}) \neq 0$.
Assume that all analytic functions of the form $\{x p + y r: (x,y) \neq 0 ; x,y \in R \}$
have only real roots. Then $\frac{G(\lambda_{k})}{F^{\prime}(\lambda_{k})} \geq 0$.
\ethm

\prf
Let $H_{++} = \{z \in C: Im(z) > 0 \}$ be the
upper half-plane. Then $\frac{G(z)}{F(z)}$ is analytic on $H_{++}$.
Also, as we explained above, the hyperbolicity of the pair $(F,G)$ implies that $Im(\frac{G(z)}{F(z)}) \neq 0 $
if $ z \in H_{++}$. Thus, from connectivity of $H_{++}$ and continuousness of 
$Im(\frac{G(z)}{F(z)})$ on $H_{++}$,
we conclude that $Im(\frac{G(z)}{F(z)})$ has the same sign on $H_{++}$. Say wlog $Im(\frac{G(z)}{F(z)}) < 0, ~z \in H_{++}$.  In other words 
$-\frac{G(z)}{F(z)}$ is a Pick function. Therefore it has the following 
integral representation \cite{bhat}:
\beqn
\frac{G(z)}{F(z)} = a + bz + \int_{-\infty}^{\infty} \frac{1+tz}{z-t} d \mu(t),
~ z \in H_{++}
\eeqn
where $ a \in R, 0 \geq b \in R$ and $\mu$ is a
nonnegative finite measure on $R$. 
Since for all $k$ a real number $\lambda_{k}$ is a simple root of $F$, 
\begin{eqnarray}
\frac{G(\lambda_{k})}{F^{\prime}(\lambda_{k})} &=& \lim_{\epsilon \downarrow 0} (\lambda_{k} + i\epsilon -\lambda_{k})
\frac{G(\lambda_{k} + i\epsilon)}{F(\lambda_{k} + i\epsilon)} = \nonumber \\
&=& \lim_{\epsilon \downarrow 0} i\epsilon \int_{-\infty}^{\infty} \frac{1+t(\lambda_{k} + i\epsilon)}{\lambda_{k} + i\epsilon-t} d \mu(t)\;. 
\end{eqnarray}

It is easy to see that
\beqn
\lim_{\epsilon \downarrow 0} i\epsilon \int_{-\infty}^{\infty} \frac{1+t(\lambda_{k} + i\epsilon)}{\lambda_{k} + i\epsilon-t} d \mu(t) =
\lim_{\epsilon \downarrow 0} i\epsilon \int_{\lambda_{k} - \delta}^{\lambda_{k} + \delta} \frac{1+t(\lambda_{k} + i\epsilon)}
{\lambda_{k} + i\epsilon-t} d \mu(t)
\eeqn
for all $\delta > 0 $. Using the identity 
$$
\frac{1 + tz}{z-t} = \frac{1 + t^{2}}{z-t} + t ;~ t,z \in C
$$
we get that 
$$
\lim_{\epsilon \downarrow 0} i\epsilon \int_{\lambda_{k} - \delta}^{\lambda_{k} + \delta} \frac{1+t(\lambda_{k} + i \epsilon)}
{\lambda_{k} + i \epsilon -t} d \mu(t) =
\lim_{\epsilon \downarrow 0} i \epsilon \int_{\lambda_{k} - \delta}^{\lambda_{k} + \delta} (\frac{1+t^{2}}{\lambda_{k} + i \epsilon -t} + t) 
d \mu(t) = 
$$

$$
=\lim_{\epsilon \downarrow 0} i\epsilon \int_{\lambda_{k} - \delta}^{\lambda_{k} + \delta} \frac{1+t^{2}}
{\lambda_{k} + i\epsilon -t}  d \mu(t)  
= \lim_{\epsilon \downarrow 0} \int_{\lambda_{k} - \delta}^{\lambda_{k} + \delta} 
\frac{(1+t^{2})(\epsilon^{2} + i\epsilon (\lambda_{k} -t)}
{\epsilon^{2} + (\lambda_{k} -t)^{2}}  d \mu(t) 
$$

As $Re(\frac{(1+t^{2})(\epsilon^{2} + i\epsilon (\lambda_{k} -t)}
{\epsilon^{2} + (\lambda_{k} -t)^{2}}) > 0$ and the last limit  exists and is real, we finally get that
$\frac{G(\lambda_{k})}{F^{\prime}(\lambda_{k})} \geq 0$. 
Notice that if $F$ and $G$ don't have common roots then $\frac{G(\lambda_{k})}{F^{\prime}(\lambda_{k})} > 0$.
\eprf

\rem It is impossible to have a hyperbolic polynomial $p(x,y,z,t) =
Det(x C_q + y C_r +z C_p -t I)$ in four variables.  Indeed, consider
three real monic polynomials $q,r,p$, all of degree $n$.  Then there
exists a real vector $(x,y,z) \neq 0$ such that $x+y+z = 01$ 
and $deg(Q) \leq n-2$ , where $Q= x q + y r + z p $.
Assume that $q$ has $n$ distinct real roots $(\lambda_{i} , 1 \leq i \leq n )$. The polynomials $q$ and $Q$ have at most $n-2$ common (real) roots :
$$
q(t) = (x-\lambda_{1})... (x-\lambda_{k}) \overline{q}(t) ,  Q(t) = (x-\lambda_{1})... (x-\lambda_{k}) \overline{Q}(t),
$$
where $deg(\overline{q}) + k =n$ and $deg(\overline{Q}) + k \leq n-2$. Therefore ,
$$
\frac{\overline{Q}(t)}{\overline{q}(t)} = \sum_{k+1 \leq i \leq n} \frac{a_{i}}{t -\lambda_{i}} ; 0 \neq a_{i} \in R , k+1 \leq i \leq n.
$$
If all $0 \neq a_{i} \in R , k+1 \leq i \leq n $ have the same sign, 
then $deg(\overline{Q}) = n-k-1$ and $deg(Q) = n-1$.
But $deg(Q) \leq n-2$ , therefore $0 \neq a_{i} \in R , k+1 \leq i \leq n $ don't have the same sign.
It follows from Theorem 3.5 that there exists $z \in H_{++}$ such that $\frac{\overline{Q}(z)}{\overline{q}(z)} = A \in R$.
This means that there exists a linear combination $ a q + Q$ with $0 \neq a \in R$ and some non-real roots. Thus
the monic polynomial of degree $n$, $P = a^{-1} ( (a+x) q + y r + z p )$! has some non-real roots and
the
determinantal polynomial $p(x,y,z,t) = Det(x C_q  + y C_r +z C_p -t I)$ is not hyperbolic in the direction $(0,0,0,1)$. 
\erem

\section{More on Alexandrov-Fenchel inequalities for mixed hyperbolic forms}
Let $p$ be an
$e$-hyperbolic polynomial of degree $n$ in $m$ variables. 
Consider $p;x_{k+1},...,x_{n}$ which are all $e$-positive. Associate with them the following polynomial of degree $k \leq n$ in
one variable
$$
\phi_{k}(t) = M_{p}(\underbrace{x+tp,...,x+tp},x_{k+1},...,x_n)
$$
Then for all $x \in R$ the roots of the polynomial $\phi$ are all real \cite{khov}. Let us call this property
$k$-hyperbolicity.
(This essentially follows from the fact that if polynomial in one variable $q$ has only real roots then its
derivative $q^{\prime}$ also has only real roots. And the latter fact is a particular case of the fact that the (complex)
roots of $q^{\prime}$ belong to the convex hull of the roots of $q$.) 

The Alexandrov-Fenchel inequalities for mixed hyperbolic forms are directly derived from 2-hyperbolicity: the
discriminant of $\phi_{2}$ is nonnegative. It is clear (for instance, from \cite{gur}) that the Alexandrov-Fenchel inequalities
(2-hyperbolicity) are not sufficient to answer Question 1, i.e. 
whether or not
$$
CO(supp(P_{x_1,..,x_n}) \cap I(n,n) = supp(Q) ?
$$
One possibility would be to use $k$-hyperbolicity for all $k \leq n$, which gives a lot of other inequalities \cite{khov}
expressed in terms of Hurwitz determinants of $\phi_{k}$ and $\phi_{k}^{\prime}$. This also might be a way to
settle Question 2 (van der Waerden conjecture for mixed hyperbolic forms). 

We will finish this section with very direct proof of concavity of $ln(p(x))$ on the positive cone $C_{e}$.

Let $x,y \in C_{e}$ and $0 < a < 1$. Then
\begin{eqnarray}
n!p(a x + (1-a)y) = \frac{1}{n!} M_{p}(a x + (1-a)y,a x + (1-a)y,...,a x + (1-a)y) \nonumber \\
= \sum_{0 \leq i \leq n}\frac{n!}{i! (n-1)!} M(i) a^{i}
(1-a)^{n-i},
\end{eqnarray}
where $M(i)$ is a  mixed hyperbolic form $M_{p}({\bf X}_{i}) $, the $n$-tuple ${\bf X}$ contains $i$ copies of $x$ and $n-i$ copies of $y$.
It follows from the Alexandrov-Fenchel inequalities that if $x,y \in C_{e}$ then \\
$0 < M(i) \geq \sqrt{M(i-1) M(i+1)}, 1 \leq i \leq n-1$. Thus $ M(i) \geq M(0)^{\frac{i}{n} } M(n)^{\frac{n-i}{n} }$,
 which gives the following inequality
\beqn
M_{p}(a x + (1-a)y,a x + (1-a)y,...,a x + (1-a)y) \geq (a M(0)^{\frac{1}{n} } + (1-a) M(n)^{\frac{1}{n} })^{n}\;.
\eeqn

Using the concavity of $ln(x), x > 0$ we get that
\beqn
ln(M_{p}(a x + (1-a)y,a x + (1-a)y,...,a x + (1-a)y)) \geq a ln(M(0)) + (1-a) ln(M(n))\;.
\eeqn
But $M_{p}(a x + (1-a)y,a x + (1-a)y,...,a x + (1-a)y) = n! p((a x + (1-a)y),  M(0) = n! p(x), M(n) = n! p(y)$, so
$$
ln(p((a x + (1-a)y)) \geq a ln(p(x)) + (1-a) ln(p(y)).
$$

\section{Acknowledgments}
First, thanks to the 
Internet, Google and Pablo Parrilo: I came across 
hyperbolic polynomials mainly because of
looking for Pablo's e-mail address. 
And Google provided... the rest.  For instance, I got the very stimulating
preprint \cite{shapiro} as a result of the
Google query "+arnold+hyperbolic+polynomial".

It is my pleasure to thank Adrian Lewis for numerous as e-mail as well phone communications.
Thanks to Leiba Rodman for discussions on the subject of {\bf Section 3}. 

I would like to acknowledge the fantastic library of Los Alamos National 
Laboratory: all references I needed were there.

\end{document}